\documentstyle[twoside]{article}
\def \C{{\Bbb C}}
\def \H{{\Bbb H}}
\def \CP{{\Bbb C}{\Bbb P}}
\def \HP{{\Bbb H}{\Bbb P}}

\def \R{{\Bbb R}}
\def \s{{\frak s}}
\def \ss{{\sigma}}

\def \k{{\kappa}}
\def \S{{\Sigma}}
\def \HH{{\mathcal H}}
\def \w{{\wedge}}

\def \o{{\omega}}

\def \sp{{\frak s}{\frak p}}
\def \so{{\frak s}{\frak o}}
\def \gl{{\frak g}{\frak l}}

\def \fd{{\bullet}}

\def \i{{\bf i}}
\def \j{{\bf j}}
\def \k{{\bf k}}
\def \proof{{\noindent{\it Proof.\ \ }}}
\def \ex{{\noindent{\it Example. \ \ }}}
\newtheorem{Th}{THEOREM}[section]
\newtheorem{definition}[Th]{DEFINITION}
\newtheorem{prop}[Th]{PROPOSITION}

\date{\it 
Department of Mathematics,
University of Arizona,
Tucson, AZ 85721-0089. 
foth@math.arizona.edu }

\title{Tetraplectic structures, tri-momentum maps, and quaternionic 
flag manifolds.}
\author{Philip Foth}
\begin{document}
\maketitle
\input amssym.def
\markboth{Philip Foth}{Tetraplectic structures, tri-momentum maps, and $\H$-flags.}
\begin{abstract} The purpose of this note is to define tri-momentum 
maps for certain manifolds with an $Sp(1)^n$-action. We exhibit many interesting
examples of such spaces using quaternions. We show how these maps can
be used to reduce such manifolds to ones with fewer symmetries. 
The images of such maps for quaternionic flag manifolds, which are defined 
using the Dieudonn\'e determinant, resemble the polytopes from the complex 
case.   
\end{abstract}

\

\noindent{\it AMS subj. class.}: \ \ primary {58F05}, secondary {53C30}.


\

\section{Introduction.}
\setcounter{equation}{0}

A closed form on a manifold can be used to gather information about its
geometry \cite{Fern}. In particular, when a manifold
is equipped with a closed non-degenerate two-form $\o$, the manifold can be
effectively studied by the methods of symplectic geometry. If, in
addition, a Lie group acts on the manifold in a hamiltonian fashion,
the Marsden-Weinstein reduction allows one to reduce the system to another 
one with fewer degrees of freedom.  

Higher order closed differential forms with properties similar to the 
symplectic structures (e.g. zero characteristic distribution)
are called multisymplectic forms and have important applications to
field theories, like in Tulczyjew \cite{Tul}, 
Marsden, Patrick, and Shkoller \cite{Marsd}, 
Cantrijn, Ibort, and de Le\'on \cite{Cant}, and other works. 
 
In this paper we concentrate on $4$-forms, and our main spaces of interest
are quaternionic vector spaces and quaternionic flag manifolds. The main 
reason is that these spaces carry natural interesting group actions, and 
appear quite naturally in many different instances. For example, quaternionic 
flag manifolds can be realized as $Sp(n)$-orbits on the 
space of quaternionic hermitian matrices. We show that these spaces carry
natural closed non-degenerate invariant $4$-forms.  

Our results are related to Yang-Mills theory. On quaternionic K\"ahler 
manifolds, Taniguchi \cite{Tanig} has
established explicit $C^0$ neighbourhoods of the minimal Yang-Mills fields,
which contain no other Yang Mills fields up to gauge equivalence.  
The quaternionic Yang-Mills connections $\nabla$ are those that satisfy 
$d_{\nabla}(F^{\nabla}\wedge \psi)=0$, where $\psi$ is a $4$-form as above. 

In general, if $X$ is an oriented $4m$-dimensional manifold equipped with a 
closed, non-degenerate $4$-form $\psi$, 
then we call $(X, \psi)$ a {\it tetraplectic} manifold. 
If, in addition, $X$ is equipped with a $Sp(1)^n$ action satisfying certain 
properties (see Section 3), we can define a {\it tri-momentum map} 
from $X$ to $(\wedge^3\s^*)^n$,
where $\s=\sp_1$. Under certain conditions, we can reduce
the original manifold $X$ to another manifold, which also possesses a 
tetraplectic structure. This procedure is quite different from the 
hyper-K\"ahler reduction of Hitchin {\it et al.} \cite{Hitch} and 
the quaternionic reduction defined by Galicki and Lawson \cite{GalLaw}, 
since our target space is different.  

We show how our procedure can be applied to the case when $X$ is a
quaternionic vector space or a (full or partial) quaternionic flag manifold. 
In particular, in the latter case, the images of the tri-momentum maps resemble
the momentum map polytops for the torus actions on the complex flag manifolds. 
The coordinate expressions for the tri-momentum maps for the quaternionic flag
manifolds can be obtained using the Dieudonn\'e determinant \cite{Dieu}. 

In the last section of the paper we discuss some related developments.

\

\noindent{\bf Acknowledgements.} I am grateful to Reyer Sjamaar for pointing
out an inaccuracy in a previous version and insightful comments. 
I am also thankful to the referee for useful comments. 
I acknowledge partial support from NSF grant DMS-0072520.

\section{Tetraplectic structures.}

We start with a definition.

\begin{definition}
Let $X$ be a real manifold of dimension $4m$ and let $\psi$ be a four-form on $X$
satisfying the following two conditions:

\noindent 1. The form $\psi$ is closed: $d\psi=0$. 

\noindent 2. The $4m$ form $\psi^m$ is the volume form on $X$. 

If these two conditions are satisfied, we call $\psi$ a \emph{tetraplectic
structure} on $X$, and $(X, \psi)$ a \emph{tetraplectic manifold}. 
\end{definition}

One of the properties of $\psi$, that is an immediate 
consequence of the definition, is that the induced maps 
$$
{\tilde{\psi}}_x: \ \ T_xX\to \w^3 T^*_xX,  
$$ 
$$
v\mapsto i_v\psi_x
$$
have trivial kernels. One finds a large class of examples of 
tetraplectic structures given by the symplectic
manifolds. If $(X, \o)$ is a $4m$-dimensional symplectic manifold, then $\psi=\o\w\o$
is a tetraplectic structure
 on $X$. However, this class of manifolds will be of little interest to us, 
since such manifolds can be effectively treated by the methods of symplectic geometry. 
More interesting examples that we have in mind include quaternionic vector spaces $\H^m$, 
full and partial quaternionic flag manifolds, and manifolds with 
$(\H^*)^n$-action. Many of the manifolds in these examples do not allow
symplectic structures.  

\

\ex
The first and basic example that we have in mind is the space $X=\H^m$. If we
identify this space with $\R^{4m}$ in the usual way, the standard 
tetraplectic form $\psi$ is defined by: 
$$
\psi=\sum_{i=1}^n dx_{4i-3}\wedge dx_{4i-2}\w dx_{4i-1}\w dx_{4i},
$$
where $x_1, ..., x_{4m}$ is the coordinate system on $\R^{4m}$. The usual 
identification between $\H^m$ with coordinates $(q_1, .., q_m)$ and $\R^{4m}$ 
is given by:
$$
q_i=x_{4i-3}+\i x_{4i-2} + \j x_{4i-1} +\k x_{4i}.
$$
We note that the form $\psi$ is not the square of a symplectic form on $\R^{4m}$
for $m >1$, because the square of a symplectic form would induce isomorphisms
$\wedge^2 T_xX \simeq \wedge^2 T^*_xX$, obtained by contraction. However, in
our case, when $m>1$, one can easily see that there will be a non-trivial
kernel at any point. 
Therefore, a naive attempt to obtain local Darboux coordinates for
every tetraplectic structure fails. Later we will discuss
a certain condition on $\psi$ which will help to get a canonical local form. 
We would like to mention in this regard that in \cite{Cant}, the authors 
describe certain canonical models of multisymplectic structures. 

\

\ex
The first important compact example of a tetraplectic structure
is given by the $4$-sphere $S^4$, which can also be
viewed as the quaternionic projective line. The tetraplectic form 
 $\psi$ on $S^4$ is just a volume form.
Under the identification with $\HP^1\simeq Sp(2)/(Sp(1)\times Sp(1))$, we can 
choose an $Sp(2)$-invariant volume form. For example, if we represent 
$S^4$ as $\H$ plus the North pole, such a form would be given by
$$
\psi = {|q|^3\over (1+|q|^4)^2}\ d|q|\ \Omega,
$$
where we identify $\H^*$ with $\R_+ \times Sp(1)$, and let $|q|$ be the
absolute value of $q\in \H$ and $\Omega$ an invariant volume 3-form on 
$Sp(1)$. (We will always view $Sp(1)$ as the group of unit length
quaternions.) Notice that $S^4$ allows neither a complex nor a  
symplectic structure. We also obtain 
a natural $Sp(1)\times Sp(1)$ action on $\HP^1$ 
coming from the natural $Sp(1)\times Sp(1)$ action on $\H^2$.   

\

\ex Let us recall the large class of quaternionic K\"ahler manifolds, 
which are $4m$-dimensional Riemannian manifolds with the holonomy 
group a subgroup of $Sp(n)Sp(1)$. It was shown by Kraines in \cite{Kra}
that all quaternionic manifolds are tetraplectic. However, this class 
does not exhaust our interest, because only quaternionic projective 
spaces are quaternionic K\"ahler, and not grassmannians or general 
flag manifolds. An interested reader should consult a beautiful
survey by Salamon \cite{Salam} and references therein. 

\

\ex A particularly large class of 4-dimensional tetraplectic
manifolds is given by the Kulkarni 4-folds \cite{Zucch}, 
\cite{Kulk}, which come 
naturally endowed with a canonical conformal class of locally conformally 
flat metrics. One can view these 4-folds as quaternionic analogues 
of Riemann surfaces. These manifolds with their volume forms 
play an important role in the 4-dimensional conformal field theory.

\begin{definition} A \emph{spheroid} $\S^n$ 
is the $n$-fold product of $Sp(1)\simeq S^3$, viewed as a Lie group. 
\end{definition}

The Lie algebra $\ss_n$ of $\S^n$ is the direct sum of $n$ copies of 
$\s=\sp_1$ - the maximal compact subalgebra of $\gl(1, \H)$. 
The Lie algebra $\s$ can be
identified with $\so(3)$ and with $\R^3$, where the Lie bracket is the
cross-product of two vectors (the vector product). 

There are natural actions of spheroids on $\H^n$ and other interesting spaces. 
This will be our main motivation for the next section. 

\

\section{Tri-momentum maps and reduction.} 

Let $X$ be a $4m$-manifold equipped with a tetraplectic structure given by
a $4$-form $\psi$. Let a spheroid $\S^n$ act on $X$ preserving the form 
$\psi$ (by tetraplectomorphisms). The stabilizer of a point $x\in X$ 
is not necessarily a spheroid. For example, if 
one considers the product $S^2\times\R^2$, and the action of $\S^1$ via
$SO(3)$ on the first factor, then there exists a volume form, which is not
changed by this action, and a stabilizer of a point is a circle. All our
further examples will be such that the stabilizer of 
are actually spheroids, and we will tacitly bear in mind this assumption for
the general discussion as well.  

If $\S$ acts on $(X, \psi)$ as above, then we have a canonical map 
$$\ss\to \Gamma(X, TX)$$ sending an element $Z$ of $\ss$
to a vector field ${\tilde Z}$ on $X$. Then we also have the map $\ss\to A^3(X)$ given by 
$Z \mapsto i_{\tilde Z}\psi$. If a three-form given by $i_Y\psi$, 
where $Y$ is a vector field on $X$, is closed, then we 
call $Y$ a \emph{locally hamiltonian vector field}. If, in addition, 
$i_Y\psi$ is exact, then we call $Y$ simply a \emph{hamiltonian vector field}. 
Here one can speculate that the group $H^3(X, \R)$ 
can be viewed as a certain topological obstruction.  

Consider the 4-vector field $\xi$ on $X$ uniquely defined by 
$i_{\xi}\psi^{m} = \psi^{m-1}$. This 4-vector field defines a quaternary
operation $\{ \cdot, \cdot, \cdot, \cdot\}$ on $C^{\infty}(X)$ 
in a standard fashion. If the Schouten 
bracket of $\xi$ with itself happens to vanish (which is the case for all 
our applications), then we get a generalized Poisson algebra structure on 
$C^{\infty}(X)$ in the terminology of \cite{Perel}. In the same source, 
as well as in \cite{Iban}, the authors consider triples of functions
$f_1$, $f_2$, $f_3\in C^{\infty}(X)$ and the
corresponding hamiltonian vector fields given by 
$$Y_{f_1,f_2,f_3}=i(df_1\wedge df_2\wedge df_3)\xi.$$ 
Then the corresponding evolution equation for any $g\in C^{\infty}(X)$
is given by $$ \dot{g} = \{ f_1, f_2, f_3, g \}.$$ 

We say that $\S$ acts on $X$ in a (generalized) 
hamiltonian way if each of the generating
vector fields for the action is hamiltonian. 
The dual vector space to the Lie algebra $\ss_n$ is isomorphic to the product of 
$n$ copies of $\s^*$. The action of $\S^n$ on $X$ induces a map 
$(\wedge^3 \s)^n\to \Gamma(X, \wedge^3TX)$ by taking the third exterior power of
the the above morphism for each component $\ss_1$ and adding these up.

\begin{definition} Let $\S^n$ act on $(X, \psi)$ in a generalized 
hamiltonian way. A \emph{tri-momentum map} $\mu$ is a map 
$$
\mu:  X\to (\wedge^3\s^*)^n\simeq \R^n
$$ 
satisfying the following conditions.

\

\noindent 1. $\mu$ is $\S^n$-invariant: $\mu(a\cdot x)=\mu(x)$, for $a\in \S^n$. 

\

\noindent 2. For any $\delta\in (\wedge^3 \s)^n$ we have
$$
d(\mu(x), \delta) = i_{\tilde{\delta}}\psi,
$$
where $x\in X$, and ${\tilde{\delta}}$ is the tri-vector field on $X$ induced by
$\delta$. 

\

\noindent 3. For any $x\in X$, such that $\mu(x)$ is regular, 
$\emph{Ker}\ T_x\mu = (\wedge^3\ss .x)^{\perp}$ with respect to $\psi_x$.
  
\end{definition}   

Notice that the first statement in the above definition is equivalent to saying 
that $\mu$ is 
$\Sigma^n$-equivariant, because the co-adjoint action of $\S^1$ on $\s^*$ induces
the trivial action of $\S^1$ on $\wedge^3\s^*$. 

We will always identify $(\wedge^3\s^*)^n$ with $\R^n$ unless it leads to 
confusion. The following is an example of a tri-momentum map. Other examples
will be treated later on. 

\

\ex
Let $X=(\H^n,\psi)$ as in Section 2 with the standard spheroid action. 
The tri-momentum map $\H^n\to \R^n$ is given by 
$$
(q_1, ..., q_n)\to (|q_1|^4, ..., |q_n|^4).
$$ 
The level sets for this tri-momentum map are isomorphic to the 
products of $3$-spheres.  

\

\ex Let us take $X=\H^2$ and the diagonal action of 
$a\in \S^1$ on $(q_1, q_2)\in \H^2$
given by $(a\cdot q_1, a\cdot q_2)$. Then the tri-momentum map 
$\H^2\to \R$ is given by $(q_1, q_2)\to |q_1|^4+|q_2|^4$. 
The regular level sets for
this tri-momentum map are isomorphic to $7$-spheres. 

\

Now we would like to define the procedure of reduction in the general setup
of tri-momentum maps. 
Let ${\bf x}=(x_1, ...,x_n)\in \R^n$ be a 
regular point of a tri-momentum map $\mu : (X, \psi)\to \R^n$ as above. The level set 
$Z_{\bf x}:=\mu^{-1}({\bf x})$ is smooth and $\S^n$-invariant. The stabilizers of the 
points in $Z_{\bf x}\subset X$ form a group bundle over it, which we assume
to be smooth. Then the 
reduced space $Y_{\bf x}:= Z_{\bf x}/\S^n$ is well defined and is a smooth manifold. 
Let us also assume that $\psi$ is \emph{horizontal} on $Z_{\bf x}$, meaning that for any 
$\beta\in \sigma_n$, and the corresponding vector field $\tilde{\beta}$ on 
$Z_{\bf x}$, one has $i_{\tilde{\beta}}\psi_{|Z_{\bf x}}=0$.  
(One can easily see that {\it a priori}, 
the 4-forms on the level sets need not necessarily be horizontal.)
By methods similar to those used for the symplectic reduction \cite{MW},
we prove the following:

\begin{Th} 
Let ${\bf x}\in \R^n$ be a regular value of a tri-momentum map $\mu: X\to \R^n$. 
Assume that the stabilizers of all points in $Z_{\bf x}$ form a smooth spheroid
bundle over $Z_{\bf x}$, and that $\psi$ is horizontal on $Z_{\bf x}$. Then the
reduced space $Y_{\bf x}=X//\S^n$ corresponding to ${\bf x}$ is a smooth
manifold admitting a tetraplectic structure $\psi_{\bf x}$, which is reduced 
from $\psi$. 
\end{Th}

\proof The tetraplectic structure $\psi$ on $X$ induces a totally
anti-symmetric 4-linear form on each of the spaces $T_z Z_{\bf x}/T_z(\S.z)$. 
This form is well defined due to the invariance and horizontality
of the form $\psi$. Therefore, we have a global four form $\psi_{\bf x}$
on the reduced space $Y_{\bf x}$. 
Now let us show that for any $y\in Y_{\bf x}$, the induced map 
$T_yY_{\bf x}\to \wedge^3 T^*_yY_{\bf x}$ has trivial kernel. 
It is enough to work with the case of $n=1$, i.e. the $Sp(1)$-action.
Since the actions of different summands commute, and are hamiltonian, 
the reduction can be performed one step at a time. In this case, one can choose a 
non-zero element in $\wedge^3\s$, which would define a Bott-Morse function $f(z)$ 
on the manifold $X$ satisfying $df=\alpha$, where $\alpha$ is the one-form,
obtained by contracting the generating tri-vector field on $X$ for the $\S^1$
action with $\psi$. According to our definition, $\alpha$ is the generating
one-form for the codimension one foliation determined by $f$ (which is 
regular, locally near the regular level sets). Therefore, we can represent 
the volume form $\psi^m$ as the product $df\wedge \beta\wedge \Omega$, 
where $\beta$ is an invariant three-form, which pairs non-trivially with the
fundamental three-vector field, and $\Omega$ is an invariant $(4m-4)$-form, 
which reduces to $Y_{\bf x}$ and is the highest exterior power of the 
tetraplectic form $\psi_{\bf x}$. $\bigcirc$

\

\ex
We leave the majority of examples for the subsequent sections, and consider only the
two examples that we had earlier in this section. 

In the first example, when $X=\H^n$ with the standard tetraplectic form 
$\psi$ and the standard $\S^n$ action, the reduced spaces are just points.  

In a slight modification of our second example, 
let $X=\H^2$ with the standard diagonal 
$\S^1$ action, and let the $4$-form $\psi$ be given by 
$$
\psi = d (|q_1|^4 - |q_2|^4)\wedge d(\alpha_1 - \alpha_2)\wedge
d(\beta_1 - \beta_2)\wedge d(\gamma_1 - \gamma_2),
$$
where $(q_1, q_2) \in \H^2$, and $q_i$ has the absolute value $|q_i|$ and the
spherical part $(\alpha_i, \beta_i, \gamma_i)$.  
The reduced space in this example 
is isomorphic to $\HP^1\simeq S^4$, and topologically 
we have the Hopf fibration $S^3\to S^7\to S^4$. The reduced
tetraplectic structure on $\HP^1$ is just the
invariant volume form discussed in Section 2. Similarly, one can obtain an
invariant tetraplectic structure on $\HP^n$ for an arbitrary $n$.  

\

We would like to reiterate that the reduction procedure described above is 
different from the Hyper-K\"ahler reduction \cite{Hitch} 
and quaternionic reduction \cite{GalLaw}. The group that acts in our
situation is $Sp(1)^n$ and the target for the momentum map involves third
exterior powers of the Lie algebra summands. Whereas, for example in 
\cite{GalLaw}, the groups maybe different, but the momentum mapping is 
bundle valued. 

We remark that one can obtain focal sets ${\rm Foc}_{\HP^n}\CP^n$ 
(critical sets for of the normal exponential map with respect to the totally
geodesic submanifold $\CP^n\subset \HP^n$) as zero level
sets of a particular momentum map as in Ornea and Piccinni \cite{OP}. It would be
interesting to see if one can obtain new examples of Sasakian-Einstein
structures using our tri-momentum maps.

\section{Quaternionic flag manifolds.}

In this section we show that the classical constructions of (full and
partial) complex flag manifolds can be used to construct quaternionic flag 
manifolds 
using the reduction procedure that we discussed in Section 3. Moreover, we
show that the reductions of these spaces possess natural invariant 
tetraplectic structures. Let
$G=GL(n, \H)$ and $B$ be the subgroup of upper triangular matrices. 
One has a natural identification between the full flag manifold 
$F_n:= Sp(n)/\S^n$ and $G/B$ similar to the complex case. We also consider 
the partial flag manifolds $F_{i_1...i_j}:=Sp(n)/(Sp(i_1)\times \cdots \times
Sp(i_j))$, where $n=i_1+\cdots +i_j$, where for example, we have 
the quaternionic grassmannians $Gr(p, n-p)=Sp(n)/(Sp(p)\times
Sp(n-p))$ appearing as conjugacy classes in the classical compact group $Sp(n)$. 
They are the orbits of the elements ${\rm diag}(\underbrace{1, ...1}_p,
\underbrace{-1, ..., -1}_{n-p})$, when $Sp(n)$ is considered as a matrix 
subgroup of $G$. The advantage of our approach is that although the lack of 
determinants over skew fields does not allow one to use the Pl\"ucker 
determinants for the
quaternionic flags, the tri-momentums maps still exist and have certain
nice properties which we will exhibit. 

Let us consider the space $\HH_n$ of quaternionic $n\times n$ 
hermitian matrices, defined as a subspace of $n\times n$ quaternionic
matrices $\gl_n(\H)$ by the condition $A=A^*$, where $A^*$ stands 
for the transposed quaternionic conjugate matrix. The group $Sp(n)$ acts
by conjugation on $\HH_n$ and the orbits of the action are isomorphic 
to quaternionic flag manifolds. 

The cell decomposition enumerated by the Schubert symbols works over $\H$ as 
well as it does over $\C$ (Ehresmann \cite{E1}). One can also use 
an identification of $\H^n$ with $\R^{4n}$ and embed quaternionic flag 
manifolds into the real ones in order to construct a non-degenerate Morse
function on $F_{i_1, ..., i_j}$, essentially done in \cite{Park}. 

A very interesting question related to the space $\HH_n$ was discussed by
Fulton in \cite{Fult}. It turns out that the equation $A_1+\cdots A_n = C$,
where the matrices have prescribed spectra, has a solution in quaternionic
hermitian matrices if and only if it has a solution in complex hermitian
matrices. 

First of all, we notice that the grassmannian $Gr(p, n-p)$ can be realized
as follows. Consider the space $\H^{np}$ of $n\times p$ matrices with 
quaternionic entries, and let the subspace $V\subset \H^{np}$ consist
of those matrices whose rows are orthonormal with respect to the
standard pairing $$\langle {\bf \xi}, {\bf \eta} \rangle = \sum_{i=1}^n
\xi_i {\bar{\eta}}_i, $$ where the bar stands for the quaternionic
conjugation. The group $Sp(p)$ acts on such bases preserving 
$V$, on which it acts freely. The quotient space is isomorphic to the 
quaternionic grassmannian $Gr(p,n-p)$. 
We claim that there exists a tetraplectic form $\psi$ on $\H^{np}$, that
can be pulled back to $V$. The resulting 4-form on $V$ will 
be preserved by the action of $Sp(p)$ and horizontal and
thus will descend to the quotient, $Gr(p, n-p)$. 
This would endow the irreducible symmetric space $Gr(p, n-p)$ with a 
tetraplectic structure. 

More generally, we can extend the construction  
of Guillemin and Sternberg in \cite{GS1} for the complex flags.
Since the full flag projects to all partial flags, and this projection 
behaves well with the respect to the group action, a 
tetraplectic structure on the full flag manifold, $F_n$, 
would push down to the partial flags. 
 
Let $K=Sp(n)$ and let $\S^n$ be its maximal spheroid as defined in Section 2. 
We have a principal fibration $K\to F_n$ with fiber $\S^n$. 
We need the following projection map 
$$ \mu_2: K\times \R^n \to \R^n,$$
$$ k \times x \mapsto x.$$ 

Here we identified for convenience
$\R^n$ with $(\wedge^3 \s^*)^n$ (see Section 3). The map $\mu_2$ is
the projection to the second factor, where we use
$\ss_n=\underbrace{\s\oplus\cdots \oplus \s}_n$ to take third exterior 
powers of individual summands. The principal bundle $K\to F_n$ admits an
invariant bundle-valued three form, $\delta$, from which we obtain a 
pre-tetraplectic 4-form 
$$
\psi_1 = d\langle \delta, \mu_2\rangle 
$$ on $K\times \R^n$ (cf. minimal coupling in \cite{Stern}). 

Following the strategy of \cite{GS1} we can show that the form $\psi_1$ is
actually tetraplectic, when we restrict to the proper subspace
$\R^n_0$ 
of $\R^n$ (using the above identification we can actually let  
$\R^n_0=(\R_+)^n$).  We call this restriction $\psi$.
Now the map $\mu_2$ which we
call $\mu$ after restricting it to $K\times \R^n_0$, has all the 
properties of a tri-momentum map from Section 3. For a generic
${\bf \xi}\in(\R_+)^n\subset \R^n$, it is clear that $\S^n$ stabilizes 
${\bf \xi}$. Therefore we obtain the following result: 

\begin{prop}
The action of $\S^n$ on $K\times \R^n_0$,  
where $\S^n$ acts trivially on the second factor, has the tri-momentum 
map $\mu$.  The  reduced space is isomorphic to 
$F_n$, the full quaternionic flag manifold. The reduced tetraplectic
form on $F_n$ so obtained is $K$-invariant. 
\end{prop}

Let us outline the relationship of the $4$-form $\psi$ on the quaternionic
flag manifolds that we have obtained with invariant symplectic structures on other 
$K$-homogeneous spaces. Let $T$ be a maximal torus in $K$ contained in $\S^n$ 
and let $T_i\simeq S^1$ be a maximal torus in the $i$-th component of 
$\S^n$, so that $T=T_1\times\cdots\times T_n$. We have the following fibration
$$ \prod_{i=1}^n(\S^1/T_i)\to Sp(n)/T \to F_n.$$
Each factor $\S^1/T_i$ is isomorphic to $\CP^1$ and carries an
$\S^1$-invariant symplectic form $\o_i$, while the space $Sp(n)/T$ is the classical 
flag manifold isomorphic to $Sp(n, \C)$ modulo its Borel subgroup, and thus carries
a $K$-invariant symplectic form $\o$. (Actually, this form can be obtained,
once one identifies $K/T$ with a co-adjoint orbit and uses the 
Kirillov-Kostant-Souriau structure on the latter.) Trivially, by choosing a fiber, 
we can assume that all the $\o_i$ are the pull-backs of the form $\o$. 
The spectral sequence for 
this fiber bundle clearly shows that, cohomologically, one can choose such a
$K$-invariant tetraplectic structure
 $\psi$ on $F_n$ that it will correspond to the cohomology 
class of $\o\wedge\o$. Moreover, due to the $K$-invariance of the
aforementioned, this correspondence can be traced on the level of forms. 

At this point, we would like to construct canonical $4$-forms on
all the orbits of $Sp(n)$ action on $\HH_n$. These orbits, as we 
mentioned earlier, are isomorphic to quaternionic flag manifolds. 
These forms have similar origins and properties to the 
KKS symplectic forms on the coadjoint orbits of the
group $U(n)$. First of all, let us define a four-commutator of square
matrices 
$$   
[A_1, A_2, A_3, A_4] = \sum_{\tau\in S_4} \emph{sign}(\tau)
A_{\tau(1)}A_{\tau(2)}A_{\tau(3)}A_{\tau(4)}.
$$
This four-commutator has the following property with respect to the usual
commutator:
\begin{equation}
\sum_{1\le i < j\le 5} (-1)^{i+j} [[A_i, A_j], A_1, ..., {\hat A_i}, ..., 
{\hat A_j}, ..., A_5]=0.
\label{eq7}
\end{equation}

We also notice that the four-commutator of four quaternionic hermitian 
matrices is again such, so we have an operation $\HH_n^{\otimes 4}\to \HH_n$.
Let us also recall the non-degenerate pairing $\HH_n\times \HH_n \to \R$
given by the real part of the trace of the product: $(A, B)\to Re\ Tr(AB)$. 
This pairing is invariant with respect to $K=Sp(n)$ action. This allows  
us to identify the tangent and co-tangent space to any element 
$y\in \HH_n$ with $\HH_n$. The 4-vector field $\kappa$ on $\HH_n$, defined
via the above 4-vector field is parallel to the orbits. The corresponding
4-form $\psi$ on an orbit ${\cal O}$, whose value at 
$y\in {\cal O}\subset \HH_n$ is given by 
$$\psi_y(A_1, A_2, A_3, A_4)=Re\ Tr(y[A_1, A_2, A_3, A_4])$$
has the following properties:

\begin{prop} The 4-form $\psi$ is non-degenerate, closed, and 
$Sp(n)$-invariant. 
\end{prop} 

\proof The invariance is a direct consequence of the fact that the
real part of the trace of the product of two quaternionic matrices is
conjugation invariant. Non-degeneracy is easy to check at one point of
the orbit, namely the diagonal matrix. Then one can use the invariance 
to show non-degeneracy on the whole orbit. 
To show that $\psi$ is closed, we will follow discussion on p.229 of 
\cite{Kir}. We will identify the orbit ${\cal O}$ with $K/L$, where
$L$ stabilizes $y$, and use the fact that $K$-invariant 4-forms on 
$K/L$ correspond uniquely to $L$-invariant elements in 
$\wedge^4({\frak l}^{\perp})$, where the differential is given by
the formula 
$$
d\phi(X_1, ..., X_5) = {1\over 5}\sum_{i < j} (-1)^{i+j+1}
\phi ([X_i, X_j], ..., {\hat X_i}, ..., {\hat X_j}, ...).
$$
Therefore, the closedness immediately follows from our formula
\ref{eq7}. $\bigcirc$

\ 

Thus, we have shown that the quaternionic flag manifolds, which appear 
as $Sp(n)$ orbits in $\HH_n$, are naturally 
tetraplectic.\footnote{I was informed that Reyer Sjammar and Yi Lin 
have a different construction of tetraplectic $4$-forms on quaternionic 
flag manifolds, using natural Lie algebra valued differential forms 
and tautological vector bundles.} 

Now we will discuss some general properties of the momentum polytopes and
we will see how the classical polytopes for the Hamiltonian torus
actions on complex flag manifolds fit into the quaternionic picture. 
The group $H=(\H^*)^n$ acts on $\H^n$ and this action induces, 
in turn, an action of $H$ on the spaces such as $F_n$ and $Gr(p, n-p)$. 
The maximal spheroid $\S^n$ is always thought of as the maximal compact 
subgroup of $H$.

Recall the Dieudonn\'e determinant \cite{Dieu}
$$
D : GL(n, \H)\to \R_+,
$$
which is defined using the transformation of a matrix 
to an upper-triangular form. For example, when $n=1$, 
$D(q)=|q|$ - the usual norm. For any $A\in GL(n, \H)$ 
of the form 
$$
\left( \begin{array}{cc}
q_1 & H \\ 0 & B \end{array} \right),
$$
where $H$ is any row vector of length $(n-1)$, and $B$ is an 
$(n-1)\times (n-1)$ matrix from $GL(n-1, \H)$,
the Dieudonn\'e determinant of $A$ is given by 
$$D(A)=|q_1|\cdot D(B).$$   

Among the properties of the Dieudonn\'e determinant are many of the
usual properties of the determinant in the group $GL(n)$ over 
a commutative field. We will use the Dieudonn\'e determinant $D$ 
to construct a tri-momentum map for the quaternionic grassmannians
$Gr(p, n-p)$. 

First of all, we will state that the combinatorics of the quaternionic
(partial) flag manifolds is not really that different from that over the
field of complex numbers $\C$. Therefore, one should expect the same major 
properties for the generators and relations of the cohomology ring 
$H^{\fd}(F_n, \C)$ and all other partial flag manifolds as in the complex
case. 

Let us first treat the case of the quaternionic grassmannian $X=Gr(p, n-p)$. 
We choose a $K$-invariant tetraplectic $4$-form $\psi$ on $X$, 
which is really only defined up to multiplication by a scalar. The spheroid
$\S^n$ acts on $X$ in a tetraplectomorphic way preserving $\psi$. 
We claim that the image of the tri-momentum map is the same as in the
complex case, i.e. can be identified with the polytope $Z^n_p$ in $\R^n$ 
defined by  
$$
Z^n_p:= \{ 0\le x_i\le 1, \ \  
x_1+\cdots + x_n=p, \}
$$
where $(x_1, ..., x_n)$ are the coordinates in $\R^n$. 
One of the ways of looking at the coordinates $x_1, ..., x_n$ is that of 
viewing them as the hamiltonians for the actions 
of the summands of $\S^n$, and we claim that the $3$-vector field
$\delta_i$, determined by the $i$-th summand in 
$ \S^n=\oplus_{i=1}^n {\S^1}$, satisfies $dx_i=i_{\delta_i}\psi$.

Now the construction of the coordinate function $x_i$ on $Gr(p, n-p)$ is not
really different from the complex case. Any quaternionc $p$-plane $\Pi$
in $\H^n$ can be viewed as an $n\times p$ matrix $M_{\Pi}$ of rank $p$
with quaternionic entries.
Following \cite{GMP}, let $J$ be a subset of $\{ 1, ..., n\}$ of cardinality
$p$. By $M_{\Pi}(J)$ we understand the $p\times p$ matrix with quaternionic 
entries obtained from $M$ by keeping only those rows that are numbered by 
the elements of $J$. We further define 
$$
x_i(\Pi)={\sum_{ i\in J} D^4(M_{\Pi}(J))\over \sum_{J}D^4(M_{\Pi}(J))}.
$$
We note that these Dieudonn\'e determinants and their 
properties are of a crucial use in the quaternionic case.

\begin{prop}
The coordinates $\{ x_i\}$ give a tri-momentum map 
$$  
\mu: Gr(p, n-p)\to \R^n.
$$ 
The image of this map is the convex polytope $Z^n_p\subset \R^n$.  
The vertices of the polytope $Z^n_p$ are the points in $Gr(p, n-p)$
fixed under the $\S^n$-action. If $\Pi$ is a point in $Gr(p, n-p)$,
then the image of the closure of the orbit $\mu(\overline{\Sigma.\Pi})$ is the
convex hull of the images of the fixed points in $\overline{\Sigma.\Pi}$. 
\end{prop}

\proof The main idea is to choose an identification between $Gr(p, n-p)$ and
an orbit of $Sp(n)$ in $\HH_n$, say of the element 
${\rm diag}(\underbrace{0,...,0}_{n-p}, \underbrace{1,...,1}_{p})$.
For example, $\HP^1$ can be identified with an orbit of $Sp(2)$ 
of ${\rm diag}(0,1)$. The element of this orbit of the form  
$$
\left( \begin{array}{cc} s_1 & s_2 \\ s_3 & s_4 \end{array} \right) \cdot 
\left( \begin{array}{cc} 0 & 0 \\ 0 & 1 \end{array} \right) \cdot
\left( \begin{array}{cc} s_1 & s_2 \\ s_3 & s_4 \end{array} \right)^{-1} = 
\left( \begin{array}{cc} |s_2|^2 & s_2 {\bar s_4} \\ s_4 {\bar s_2} & |s_4|^2 
\end{array} \right), 
$$ 
where the first matrix is from $Sp(2)$, and the second from $\HH_2$,
corresponds to $[s_2 : s_4]$ $\in \HP^1$. On the other hand, a point 
$\HP^1$ can be represented as a line in $\H^2$ passing through a point   
$(q_1, q_2)$. There are certain formulas relating $q_1, q_2$ with 
$s_1, s_3$, which can be obtained by a straightforward computation.  
Analogous considerations are valid for all $Gr(p, n-p)$'s. 
One the orbit side, the momentum map will simply be given by the projection 
to the diagonal. 
$\bigcirc$

\

Similarly, one can obtain statements about 
full and all partial quaternionic flag manifolds 
that are analogous to the
complex case.

\section{Related developments.}

In this section we will merely outline the content of two subsequent papers
\cite{F2} and \cite{FL}. 
In \cite{F2} we show several important generalizations of the classical
results from symplectic geometry to the case of tetraplectic geometry. 
The first is the convexity theorem founded in \cite{At} and \cite{GS3}. 
Basically, we establish the following fact. Let $(X, \psi)$ be a tetraplectic
manifold and let $\S^n$ act on $X$ in a tetraplectomorphic
way. Let $(f_1, ..., f_n)$ be such functions on $X$ that the flow
corresponding to the generalized hamiltonian $3$-vector fields 
$\delta_1$, ..., $\delta_n$ (defined by $i_{\delta_j}\psi=df_j$) generated
a subgroup of {\it Diff}$(X)$ defined by $\S^n$. Then the image of the
map $\mu: X\to \R^n$ given by 
$$
\mu({\bf x})=(f_1({\bf x}), ..., f_n({\bf x}))
$$ is the convex hull of the images of connected components of the set
of common critical points of $f_i$'s.

Another direction that we pursued in \cite{F2} is a generalization of the
Duistermaat - Heckman theorems \cite{DH}. For the case of quaternionic 
flag manifolds and certain other compact manifolds with $(\H^*)^n$ action,
this formula would help to recover the structure of the cohomology
ring of the  reduced spaces
from the combinatorics of the fixed point data combined with 
a generalization of Duistermaat-Heckman by methods similar to 
Guillemin-Sternberg \cite{GS1}. 

In \cite{FL} we work with generalized Poisson structures (GPS) of rank 4 as defined 
by de Asc\'arraga, Perelomov, and P\'erez Bueno in \cite{Perel}. We show that
many familiar manifolds have natural GPS. In particular, we show that the full
quaternionic flag manifolds $F_n$ (as well as all the partial ones) 
have interesting natural GPS. In particular, we give a Lie theoretic
construction of the Bruhat 4-vector field on $F_n$, which is an example of GPS.
Recall that the classical Bruhat Poisson structures 
on complex flag manifolds that were first
introduced by Soibelman in \cite{Soi} and independently by Lu and Weinstein
in \cite{LW}. One of their main properties is that
that the symplectic leaf decomposition yields exactly the
Bruhat cells. Moreover, Evens and Lu in \cite{EL} showed that the Kostant
harmonic forms from \cite{Kost} have Poisson harmonic nature with respect 
to the Bruhat Poisson structure. 
In \cite{FL} we show that there is a basis in cohomology of 
$H^\fd(F_n)$ dual to the natural Bruhat cell decomposition for quaternionic
flag manifolds $F_n$ that is represented by forms with properties similar
to Kostant harmonic forms. In particular, those lead to harmonic forms with 
respect to our new Bruhat 4-vector field. It is also quite natural to consider the 
equivariant cohomology with the respect to the natural spheroid action on
$F_n$. 

We also plan to study analogues of other interesting facts from the 
complex geometry, which can be adapted to the quaternionic case. 
Examples include the Gelfand-MacPherson correspondence between GIT
and symplectic quotients of grassmannians and products of projective 
spaces, the Gelfand-Tsetlin coordinates on the space of hermitian matrices, 
moduli spaces of quaternionic vector bundles, and others. 

\

\noindent{\it Manifolds with $(\H^*)^n$-action. \ }
One could be tempted to use the theory of the toric manifolds
to study the manifolds with an $(\H^*)^n$ action with a dense open orbit. 
In particular, one can start with a convex polytope in $\R^n$
and try to construct a $4n$-dimensional manifold with an
$(\H^*)^n$ action that has a dense open orbit such that the
tri-momentum map for the corresponding $\S^n$-action is that
convex polytope. We do not know if this is possible in general 
except the simplest situation, when the polytope is the standard simplex
in $\R^n$. In this case, the corresponding manifold is $\HP^n$. 
A naive application of the  usual reduction method of constructing 
such manifolds fails in general due to the non-commutative nature of 
$\S^1$. However, there are examples of classes of manifolds
on which $(\H^*)^n$ acts with a dense open orbit and we believe that 
one can classify those by methods similar to \cite{Dan} and \cite{Del}. 
Scott \cite{Scott} develops a theory of quaternionic toric manifolds, 
using topological methods. In general, only a single copy of $Sp(1)$ acts
on those. If we only assume that $\H^*$ acts on a manifold $X$, then in certain 
cases one can obtain a cell decomposition of $X$ similar to the Bialynicki-Birula 
decomposition in the complex case. Examples of such spaces are given by 
the quaternionic flag manifolds.


\thebibliography{123}

\bibitem{At}{M. Atiyah. Convexity and commuting hamiltonians.
{\it Bull. London Math. Soc.}, {\bf 14}:1-15, 1982.}

\bibitem{Au}{M. Audin. The topology of torus actions on symplectic manifolds.
{\it Birkh\"auser Verlag}, 1993.}

\bibitem{Cant}{F. Cantrijn, A. Ibort, and M. de Le\'on. On the geometry of 
multisymplectic manifolds. {\it J. Austral. Math. Soc. Ser. A}, {\bf 66}:
303-330, 1999.}

\bibitem{Dan}{V. I. Danilov. The geometry of toric varieties. 
{\it Russian Math. Surveys}, {\bf 33}: 97-154, 1978.}

\bibitem{Perel}{J. A. de Asc\'arraga, A. M. Perelomov, and 
J. C. P\'erez Bueno. The Schouten-Nijenhuis bracket, cohomology
and generalized Poisson structures. {\it J. Phys. A}, {\bf 29}: 
7993-8009, 1996.}

\bibitem{Del}{T. Delzant. Hamiltoniens p\'eriodiques et image convexe de
l'application moment. In French. {\it Bull. Soc. Math. France}, 
{\bf 116}: 315-339, 1988.}

\bibitem{Dieu}{J. Dieudonn\'e. Les d\'eterminants sur un corps non
commutatif. In French. {\it Bull. Soc. Math. France}, 
{\bf 71}: 27-45, 1943.}

\bibitem{DH}{J.J. Duistermaat and G. Heckman. On the variation of cohomology 
of the symplectic form of the reduced phase space. {\it Invent. Math.}, 
{\bf 69}; 259-269, 1982.}  

\bibitem{E1}{C. Ehresmann. Sur la topologie de certains espaces homog\'enes. In
French. {\it Annals of Math.}, {\bf 35}: 396-443, 1935.}

\bibitem{EL}{S. Evens and J.-H. Lu. Poisson harmonic forms, Kostant harmonic
forms, and the $S^1$-equivariant cohomology of $K/T$. {\it Advances in Math.}, 
{\bf 142}: 171-220, 1999.} 

\bibitem{Fern}{M. Fern\'andez, R. Ib\'a\~nez, and M. de Le\'on. The geometry of
a closed form. {\it Homotopy and geometry. Warsaw, 1997}, Banach
Center Publ., {\bf 45}: 155-167, 1997.} 

\bibitem{FL}{P. Foth and F. Leitner. Geometry of four-vector fields on 
quaternionic flag manifolds. {\it math.DG/0104104.}}

\bibitem{F2}{P. Foth. Cohomology of tetraplectic quotients.
{\it In preparation.}} 

\bibitem{Fult}{W. Fulton. Eigenvalues, invariant factors, heighest weights, and
Schubert calculus. {\it Bull. Amer. Math. Soc.}, {\bf 37}: 209-249, 2000.}

\bibitem{GalLaw}{K. Galicki and H. B. Lawson, Jr. Quaternionic reduction 
and quaternionic orbifolds. {\it Math. Ann.}, {\bf 282}: 1-21, 1988.}

\bibitem{GMP}{I. M. Gelfand, R. Goresky, R. MacPherson, and V. V. Serganova. 
Combinatorial geometries, convex polyhedra, and Shubert cells. 
{\it Advances in Math.}, {\bf 63}: 301-316, 1987.}  

\bibitem{GS1}{V. Guillemin and S. Sternberg. The coefficients of the 
Duistermaat-Heckman polynomial and the cohomology ring of reduced spaces.
{\it Geometry, topology, and physics}: 202-213, 
International Press, Cambridge, MA, 1995.}

\bibitem{GS3}{V. Guillemin and S. Sternberg. Convexity properties of the
moment maps. {\it Invent. Math.}, {\bf 97}: 491-513, 1982.}

\bibitem{Hitch}{N. Hitchin, A. Karlhede, U. Lindstr\"om, and M. Ro\v{c}ek. 
Hyper-K\"ahler metrics and supersymmetry. {\it Comm. Math. Phys.}, 
{\bf 108}: 535-589, 1987.}

\bibitem{Iban}{R. Ib\'a\~nez, M. de Le\'on, J. Marrero, and D. M. de Diego. 
Dynamics of generalized Poisson and Nambu-Poisson brackets.   
{\it J. Math. Phys.}, {\bf 38}: 2332-2344, 1997.}   

\bibitem{Kir}{A. A. Kirillov. Elements of the theory of representations. 
{\it Grundlehren der Mat. Wiss.}, {\bf 220}, Springer-Verlag, 1976.}

\bibitem{Kra}{V. Y. Kraines. Topology of quaternionic manifolds. 
{\it Trans. AMS}, {\bf 122}: 357-367, 1966.} 

\bibitem{Kost}{B. Kostant. Lie algebra cohomology and generalized Shubert cells.
{\it Ann. Math.}, {\bf 77}: 72-144, 1963.}

\bibitem{Kulk}{R. S. Kulkarni. On the principle of uniformization. 
{\it J. Diff. Geom.}, {\bf 13}: 109-138, 1978.}  

\bibitem{LW}{J.-H. Lu and A. Weinstein. Poisson-Lie groups, 
dressing transformations, and Bruhat decompositions. {\it J. Diff. Geom.}, 
{\bf 31}, 501-526, 1990.}

\bibitem{Marsd}{J. Marsden, G. Patrick, and S. Shkoller. Multisymplectic
geometry, variational integrators, and nonlinear PDEs. {\it Comm. Math. Phys.},
{\bf 199}: 351-395, 1998.}

\bibitem{MW}{J. Marsden and A. Weinstein. Reduction of symplectic manifolds
with symmetry. {\it Rep. Mathematical Phys.}, {\bf 5}: 121-130, 1974.}

\bibitem{OP}{L. Ornea and P. Piccinni. On some moment maps and induced Hopf
bundles in the quaternionic projective space. {\it Internat. J. Math.}, 
{\bf 11}: 925-942, 2000.}  

\bibitem{Park}{G. Parker. Morse theory on quaternionic grassmannians. 
{\it J. Diff. Geometry}, {\bf 7}: 615-619, 1972.} 

\bibitem{Salam}{S. Salamon. Quaternion-K\"ahler geometry. In 
{\it Surveys in Differential Geometry IV: Essays on Einstein Manifolds}:
83-123, International Press, 1999.} 

\bibitem{Soi}{Y. Soibelman. 
The algebra of functions on a compact quantum group 
and its representations. {\it Leningrad J. Math.}, {\bf 2}: 161-178, 1991.}

\bibitem{Scott}{R. Scott. Quaternionic toric varieties. {\it Duke Math. J.}, 
{\bf 78}: 373-397, 1995.}

\bibitem{Stern}{S. Sternberg. Minimal coupling and symplectic mechanics
of a classical particle in the presence of a Yang-Mills field. {\it
Proc. Nat. Acd. Sci. U.S.A.}, {\bf 74}: 5253-4, 1977.} 

\bibitem{Tanig}{T. Taniguchi. Isolation phenomena for quaternionic Yang-Mills 
connections. {\it Osaka J. Math.}, {\bf 35}: 147-164, 1998.}  

\bibitem{Tul}{W. M. Tulczyjew. Homogeneous symplectic formulation of field
dynamics and the Poincar\'e-Cartan form. {\it Lect. Not. Math.}, {\bf 836}:
484-497, Springer, Berlin, 1980.}

\bibitem{Zucch}{R. Zucchini. The quaternionic geometry of four-dimensional
conformal field theory. {\it J. Geom. Phys.}, {\bf 27}: 113-153, 1998.}
\end{document}